\documentclass{amsart}
\usepackage{amsfonts,amssymb,amsmath,amsthm}
\usepackage{enumerate,latexsym}
\usepackage{amsmath, amsthm, amscd, amsfonts,url, amssymb,graphicx,graphics, color}

\newtheorem{thm}{Theorem}[section]

\newtheorem{cor}[thm]{Corollary}
\theoremstyle{definition}

\newtheorem{rem}[thm]{Remark}
\numberwithin{equation}{section}

\author{A. Shirafkan}
\address{Department of Mathematics, Faculty of Mathematics Sciences,\\
University of Mazandaran, Babolsar, Iran.}
\email{shirafkan@umz.ac.ir}
\author{M. Rafie-Rad}
\address{Department of Mathematics, Faculty of Mathematics Sciences,\\
University of Mazandaran, Babolsar, Iran.}
\email{rafie-rad@umz.ac.ir}

\begin{document}

\title[On the C-projective vector fields on Randers spaces]{On the C-projective vector fields on Randers spaces}

\begin{abstract}
A characterization of the C-projective vector fields on a Randers
spaces is presented in terms of a recently introduced non-Riemannian
quantity defined by Z. Shen and denoted by ${\bf\Xi}$; It is proved
that the quantity ${\bf\Xi}$ is invariant for C-projective vector
fields. Therefore, the dimension of the algebra of the C-projective
vector fields on an $n$-dimensional Randers space is at most
$n(n+2)$. The generalized Funk metrics on the $n$-dimensional
Euclidean unit ball $\mathbb{B}^n(1)$ are shown to be explicit
examples of the Randers metrics with a C-projective algebra of
maximum dimension $n(n+2)$. Then, it is also proved that an
$n$-dimensional Randers space has a C-projective algebra of maximum
dimension $n(n+2)$ if and only if it is locally Minkowskian or (up
to re-scaling) locally isometric to the generalized Funk metric. A
new projective invariant is also introduced.
\end{abstract}
 \maketitle
\section{Introduction}
The projective and conformal structures on a manifolds determine the
metric structure. The projective properties of a generic Finsler
space may include a sever inconvenience in comparison to the
Riemannian spaces: the Beltrami's theorem in Riemannian geometry is
no longer valid in Finsler geometry (cf. \cite{BaSh}) and this turns
the projective Finsler geometry Bizarr. In a Finsler space $(M,F)$,
the usual Ricci tensor $K_{jl}=K^i_{\ jil}$ (with respect to the
Berwald connection) is not generally symmetric with respect to the
indices $j$ and $l$. Therefore, an infinitesimal projective
transformation may convert or even preserves the anti-symmetric part
$\mathcal{R}_{jl}=\big(K_{jl}-K_{lj}\big)/2$ which is in turn a
non-Riemannian quantity. Every projective vector field on a
Riemannian background space preserves $\mathcal{R}$ trivially and it
can be realized as the so called \textit{closeness property}.
However, in general a projective vector field on a generic Finsler
space may refuse the have the closeness property. A projective
vector field is \textit{C-projective} if it preserves $\mathcal{R}$
(or in other words it has the closeness property). In the recent
work \cite{NajafiTayebi} the Lie algebra of the C-projective vector
fields has been studied and a C-projective invariant has been also
introduced. The well-known non-Riemannian
curvature {\bf H}-curvature is an invariant of the algebras of C-projective vector fields.\\
\bigskip
Given any Finsler metric $F$ on an $n$-dimensional manifold $M$,
consider the following four non-Riemannian S-curvature driven
quantities ${\bf \Xi}= \Xi_idx^i$, ${\bf E}=E_{ij}dx^i\otimes dx^j$,
${\bf H}= H_{ij}dx^i\otimes dx^j$ and ${\bf
\Sigma}=\Sigma_{ij}dx^i\otimes dx^j$ on the pullback tangent bundle
$\pi^*TM$:
\begin{eqnarray}
\Xi_i&=&{\bf S}_{.i|m}y^m-{\bf S}_{|i},\\
E_{ij}&=&\frac{1}{2}{\bf S}_{.i.j}\\
H_{ij}&=&\frac{1}{2}{\bf S}_{.i.j|m}y^m,\\
\Sigma_{ij}&=&\frac{1}{n+1}\Big({\bf S}_{.i|j}-{\bf S}_{.j|i}\Big),
\end{eqnarray}
where ${\bf S}$ denotes the S-curvature and ``." and ``$|$" denote
the vertical and horizontal covariant derivatives, respectively,
with respect to the Berwald connection. The quantity $\Xi$ was first
introduced by Shen in \cite{Shen4}. In fact, the above quantities do
not depend to the choice of connection for performing horizontal
derivatives and can be derived for the Finsler metric itself. Notice
that the following implications for Randers metrics are useful (cf.
\cite{Shen4}):
\[
{\bf \Sigma}=0\ \ \Leftrightarrow\ \ {\bf \Xi}=0\ \ \Leftrightarrow\
\ {\bf S}=(n+1)cF,\ (c\in\mathbb{R}).
\]
Here, we prove the following result:
\begin{thm}
\label{mainthm1}
 Let $(M,F=\alpha +\beta )$ be a Randers space and $V$ is a projective vector field
$V\in \chi (M)$. The following statements are equivalent:\\
(1) $V$ is C-projective,\\
(2) $\mathcal{L}_{\hat{V}}{\bf \Sigma}=0$,\\
(3) $\mathcal{L}_{\hat{V}}{\bf \Xi}=0$.
\end{thm}
Theorem \ref{mainthm1} ensures that $\Xi$ and $\Sigma$ are
C-projectively invariant quantities. On an $n$-dimensional Randers
space $(M,F=\alpha+\beta)$, denote the algebra of projective and
C-projective vector fields by $p(M,F)$ and $cp(M,F)$, respectively.
Then, $\dim\big(cp(M,F)\big)\leq n(n+2)$ as well as the projective
algebra. However, the case of maximum dimension is also interesting:
The locally Minkowski spaces are trivial examples possessing the
C-projective algebra of maximum dimension. The following Randers
metric on the Euclidean unit ball $\mathbb{B}^n(1)$ is called the
generalized Funk metric:
\begin{equation}
\label{local form} F(x,y)=\frac{\sqrt{|y|^2-(|x|^2|y|^2-\langle
x,y\rangle^2)}}{1-|x|^2}\pm\frac{\langle
x,y\rangle}{1-|x|^2}\pm\frac{\langle a,y\rangle}{1+\langle
a,x\rangle},
\end{equation}
where, $\langle,\rangle$ and $|.|$ denote the Euclidean inner
product and norm on $\mathbb{R}^n$, $y\in T_x\mathbb{R}^n,\
a\in\mathbb{R}^n, |a|<1$. The C-projective algebra takes its maximum
dimension $n(n+2)$ for the generalized Funk metrics.

\bigskip
\begin{thm}
\label{mainthm2} Let $(M,F=\alpha +\beta )$ be a Randers space of
dimension $n\geq3$. The C-projective algebra $cp(M,F)$ takes the
maximum (local) dimension $n(n+2)$, if and only if $F$ is a locally
Minkowski metric or it is (up to re-scaling) locally isometric to
the generalized Funk metric on the Euclidean unit ball
$\mathbb{B}^n(1)$

\end{thm}
The horizontal derivatives $D_{\frac{\delta}{\delta x^i}}$ with
respect to the Berwald connection $D$ is denoted by $_{|i}$. The
subscripts $_{;i}$ and $_{.i}$ stand for the partial derivations
$\frac{\partial}{\partial x^i}$ and $\frac{\partial}{\partial y^i}$
are , respectively. The complete lift of any vector field $V$ on $M$
o $TM$ is denoted by $\hat{V}$ and $\mathcal{L}_{\hat{V}}$ denotes
the Lie derivative operator with respect to $\hat{V}$. Moreover, we
deal with pure Randers metrics, i.e $\beta\neq0$.
\section{Preliminaries}\label{sectionP}
Let $M$ be a $n$-dimensional $ C^\infty$ connected manifold. The
tangent space of $M$ at $x\in M$ is denoted by $T_x M $ and the
tangent manifold of $M$ is the disjoint union of tangent spaces
$TM:=\bigsqcup _{x \in M} T_x M$. Every element of $TM$ is a pair
$(x, y)$ where $x\in M$ and $y\in T_xM$. Denote the slit tangent
manifold by $TM_0 = TM\setminus \{\bf o\}$, where, ${\bf o}$ denotes
the zero section of the tangent bundle.. The natural projection
$\pi: TM \rightarrow M$ given by $\pi (x,y):= x$ makes $TM$ a vector
bundle of rank $n$ over $M$ and $TM_0$ a fiber bundle over $M$ with
fiber type $\mathbb{R}^n\setminus\{\bf o\}$. A {\it Finsler metric}
on $M$ is a function $ F:TM \rightarrow [0,\infty )$ satisfying
following conditions: (i) $F$ is $C^\infty$ on $TM_0$,\ (ii)
$F(x,y)$ is positively 1-homogeneous $y$ and (iii) the Hessian
matrix of $F^{2}$ with entries $
g_{ij}(x,y):=\frac{1}{2}[F^2(x,y)]_{y^iy^j} $ is positively defined
on $TM_0$.  Given any Finsler metric $F$ on $M$, the pair $(M,F)$ is
called a {\it Finsler space}. Traditionally, we denote a Riemannian
metric by $\alpha=\sqrt{a_{ij}(x)y^iy^j}$. The geodesic spray $G$ is
naturally induced by $F$ on $TM_0$ given in any standard coordinate
$(x^i,y^i)$ for $TM_0$ by $\textbf{G}=y^i {{\partial} \over
{\partial x^i}}-2G^i(x,y){{\partial} \over {\partial y^i}}$, where
$G^i(x,y)$ are local functions on $TM_0$ given by
$G^i:=\frac{1}{4}g^{ih}\Big\{y^kF^2_{x^ky^h}-F^2_{x^h}\Big\}.$
Assume the following conventions:
\[
G^i_{\ j}=\frac{\partial G^i}{\partial y^i},\ \ \ G^i_{\
jk}=\frac{\partial G^i_{\ j}}{\partial y^k},\ \ \
\]
The local functions $G^i_{\ j}$ are coefficients of a connection in
the pullback bundle $\pi^*TM\longrightarrow M$ which is called the
\textit{Berwald connection} denoted by $D$. Recall that for
instance, the derivatives of a vector field $V$ and a 2-covariant
tensor $T=T_{ij}dx^i\otimes dx^j$ is given by:
\begin{eqnarray}
X^i_{\ |k}&=&\frac{\delta X^i}{\delta x^k}+X^rG^i_{rk},\nonumber\\
X_{i|k}&=&\frac{\delta X_i}{\delta x^k}-X_rG^r_{ki},\label{Cov der}\\
T_{ij|k}&=&\frac{\delta T_{ij}}{\delta
x^k}-T_{rj}G^r_{ik}-T_{ir}G^r_{kj},\nonumber
\end{eqnarray}
where, $\frac{\delta}{\delta x^k}=\frac{\partial}{\partial
x^k}-G^i_{\ k}\frac{\partial}{\partial y^i}$. The {\it
Busemann-Hausdorff volume form} $dV_F = \sigma_F(x) dx^1 \cdots
dx^n$ on any Finsler space $(M,F)$ is defined by
\[
\sigma_F(x) := \frac{\textrm{Vol} (\Bbb B^n(1))}{ \textrm{Vol}  \{
(y^i)\in \mathbb{R}^n \  | \ F  ( y^i \frac{\partial}{\partial
x^i}|_x ) < 1 \} }.
\]
The {\bf S}-curvature with respect to the Busemann-Hausdorff volume
form is denoted by ${\bf S}$ and is defined by
\begin{equation}
\label{S spray} {\bf S}=G^m_{\ m}-y^m\frac{\partial}{\partial
x^m}\ln\sigma_F.
\end{equation}
Given any Randers space $(M,F=\alpha+\beta)$, the ${\bf
S}$-curvature takes the following form:
\begin{equation}
\label{Randers} {\bf
S}=(n+1)\big\{\frac{e_{00}}{2F}-s_0-\rho_0\big\},
\end{equation}
where, $\rho=\ln(\sqrt{1-\|\beta_x\|^2_\alpha})$ and
$\rho_0=y^ic_{x^i}$. We may consider several quantities using the
{\bf S}-curvature. The $\Xi$, $\Sigma$ and $H$ curvatures are
denoted by $\Xi$ and ${\bf \Sigma}$ and ${\bf H}$ respectively, and
are defined at every point $x\in M$ by
\begin{eqnarray*}
{\bf\Xi}_y&=&\Xi_i(y)dx^i,\\
{\bf E}_y&=&E_{ij}dx^i\otimes dx^j \\
{\bf H}_y&=&H_{ij}(y)dx^i\otimes dx^j\\
{\bf \Sigma}_y&=&\Sigma_{ij}(y)dx^i\otimes dx^j,\\
\end{eqnarray*}
where $y\in T_xM-\{0\}$ and
\begin{eqnarray}
\Xi_i&:=&y^m{\bf S}_{.i|m}-{\bf S}_{|i},\label{Sc 01}\\
E_{ij}&:=&\frac{1}{2}{\bf S}_{.i.j},\label{Sc 02}\\
H_{ij}&:=&\frac{1}{2}y^m{\bf S}_{.i.j|m},\label{Sc 03}\\
\Sigma_{ij}&:=&\frac{1}{n+1}\Big\{{\bf S}_{.i|j}-{\bf
S}_{.j|i}\Big\}.\label{Sc 04}
\end{eqnarray}
Although the ${\bf S}$-curvature depends to the chosen volume form,
following the equations \eqref{S spray},\eqref{Sc 01},\eqref{Sc
02},\eqref{Sc 03} and \eqref{Sc 04}, it not hard to show that ${\bf
\Xi}$, ${\bf E}$, ${\bf \Sigma}$ and ${\bf H}$ are independent from
choosing any volume form for the {\bf S}-curvature; For example, we
may see it for ${\bf\Sigma}$ below: By the definition of the {\bf
S}-curvature ${\bf S}=\Pi-y^m\frac{\partial}{\partial
x^m}\ln\sigma_F$, we have:
\begin{eqnarray*}
(n+1)\Sigma_{ij}&=&{\bf
S}_{.i|j}-{\bf S}_{.j|i}=\frac{\delta}{\delta x^j}{\bf S}_{.i}-{\bf S}_{.k}G^k_{\ ji}-\frac{\delta}{\delta x^i}{\bf S}_{.j}+{\bf S}_{.k}G^k_{\ ij}\\
&=&\frac{\delta}{\delta x^j}{\bf S}_{.i}-\frac{\delta}{\delta
x^i}{\bf S}_{.j}\\
&=&\frac{\delta}{\delta x^j}\big(\Pi_{y^i}-\frac{\partial}{\partial
x^i}\ln\sigma_F\big)-\frac{\delta}{\delta
x^i}\big(\Pi_{y^j}-\frac{\partial}{\partial x^j}\ln\sigma_F\big)\\
&=&\big(\frac{\delta}{\delta x^j}\Pi_{y^i}-\frac{\delta}{\delta
x^i}\Pi_{y^j}\big)+\big(\frac{\partial^2}{\partial
x^ix^j}\ln\sigma_F-\frac{\partial^2}{\partial
x^jx^i}\ln\sigma_F\big)\\
&=&\Pi_{y^ix^j}-G^r_{\ j}\Pi_{y^ry^i}-\Pi_{y^jx^i}-G^r_{\
i}\Pi_{y^ry^j},
\end{eqnarray*}
which follows
\[
\Sigma_{ij}=\frac{1}{n+1}\Big\{\Pi_{y^ix^j}-\Pi_{y^ry^j}G^r_{\
i}-\Pi_{y^jx^i}+\Pi_{y^ry^i}G^r_{\ j}\Big\}
\]
and $\Pi=G^m_{\ m}$. It is not hard to show that here are fine
relations between the above four quantities given below:
\begin{eqnarray}
& &\Sigma_{ij}=-\Sigma_{ji}\label{R00}\\
& &y^i{\Sigma}_{ij}=-\frac{1}{n+1}\Xi_j\label{R01}\\
& &y^j\Xi_{j.k}=-\Xi_k\label{R02}\\
& &\Xi_{i.j}+\Xi_{j.i}=4H_{ij},\ \ (\textrm{cf.}\  [9])\label{R1}\\
& &\Xi_{i.j}-\Xi_{j.i}=2(n+1)\Sigma_{ij},\label{R2}\\
& &y^i\Sigma_{ij.k}=-\frac{2}{n+1}H_{jk}, \label{R3}\\
& &\Xi_{j.k}=2H_{jk}+(n+1)\Sigma_{jk}, \label{R4}
\end{eqnarray}
And from \eqref{R02} and \eqref{R4} the following equations result:
\begin{eqnarray}
& &-\Xi_k=y^j\Xi_{j.k}=(n+1)y^j\Sigma_{jk}.\label{R5}
\end{eqnarray}
 A Finsler space is said to be {\it of isotropic {\bf
S}-curvature} if there is a function $c=c(x)$ defined on $M$ such
that ${\bf S}=(n+1)c(x)F$. It is called a Finsler space {\it of
constant {\bf S}-curvature} once $c$ is a constant. Every Berwald
space is of vanishing {\bf S}-curvature \cite{Shen2}. The following
result proved by Z. Shen show that constancy of S-curvature and
vanishing of $\Xi$ re the same for Randers metrics:
\bigskip
\begin{thm}\label{xi S constant} (Z. Shen, \cite{Shen4})
\label{Xi=0} Let $F=\alpha+\beta$ be a Randers metric on an
$n$-dimensional manifold $M$.Then, ${\bf S}=(n+1)cF$ for some
constant $c$ if and only if ${\bf \Xi}=0$.
\end{thm}
Let $(M,\alpha)$ be a Riemannian space and $\beta=b_i(x)y^i$ be a
1-form defined on $M$ such that $\|\beta_x\|_\alpha :=\sup_{y \in
T_xM} \beta(y)/\alpha(y) <1$. The Finsler metric $F = \alpha+\beta$
is called a Randers metric on a manifold $M$. Denote the geodesic
spray coefficients of $\alpha$ and $F$ by the notations $G_\alpha^i$
and $G^i$, respectively and the Levi-Civita connection of $\alpha$
by $\widetilde{\nabla}$. Define $\widetilde{\nabla}_jb_i$ by
$(\widetilde{\nabla}_jb_i) \theta^j := db_i -b_j \theta_i^{\ j}$,
where $\theta^i :=dx^i$ and $\theta_i^{\ j} :=\tilde{\Gamma}^j_{ik}
dx^k$ denote the Levi-Civita connection forms and
$\widetilde{\nabla}$ denotes its associated covariant derivation of
$\alpha$. Recall the conventional standard notations for Randers
metrics given by $r_{ij}:={1\over
2}(\widetilde{\nabla}_jb_i+\widetilde{\nabla}_ib_j), s_{ij}:=
{1\over 2}(\widetilde{\nabla}_jb_i-\widetilde{\nabla}_ib_j), s^i_{\
j}:= a^{ih}s_{hj}, s_j:=b_i s^i_{\ j}$ and $e_{ij} := r_{ij}+ b_i
s_j + b_j s_i$. Then $G^i$  is given by
\begin{equation}
G^i = G_\alpha^i + \Big({e_{00} \over 2F} -s_0\Big)y^i+ \alpha
s^i_{\ 0},\label{Gii}
\end{equation}
where $e_{00}:= e_{ij}y^iy^j$, $s_0:=s_iy^i$, $s^i_{\ 0}:=s^i_{\ j}
y^j$ and $G^i_\alpha$ denote the geodesic coefficients of $\alpha$,
see \cite{Shen2}. It is well-known that a Randers metric $F$ is of
isotropic {\bf S}-curvature ${\bf S}=(n+1)c(x)F$ if and only if
$e_{00}=2c(x)(\alpha^2-\beta^2)$; see \cite{ShenXing}. Therefore,
for a Randers metric of isotropic {\bf S}-curvature the spray
coefficients $G^i$ are of the form
\begin{equation}
G^i = G_\alpha^i + \Big(c(x)(\alpha-\beta)-s_0\Big)y^i+\alpha s^i_{\
0},\label{Gi}
\end{equation}
Notice that, due to \eqref{Randers}the coefficients $G^i$ can be
written in terms the {\bf S}-curvature as follows:
\begin{equation}
G^i = G_\alpha^i +\Big(\frac{{\bf S}}{n+1}+\rho_0\Big)y^i+ \alpha
s^i_{\ 0},\label{GiS}
\end{equation}

The Riemann curvature tensor is defined by ${\bf R}_y= R^i_{\
k}(x,y) dx^k \otimes \frac{\partial}{\partial x^i}|_x: T_xM
\longrightarrow T_xM$ by
\[
 R^i_{\ k} := 2 \frac{\partial G^i}{\partial x^k}
- y^j \frac{\partial^2 G^i}{\partial x^j\partial y^k} + 2 G^j
\frac{\partial^2 G^i}{\partial y^j\partial y^k}-\frac{\partial
G^i}{\partial y^j} \frac{\partial G^j}{\partial y^k}.
\]
The family ${\bf R}:=\Big\{{\bf R}_{\pi z}\Big\}_{z\in TM_0}$ is
called the Riemann curvature \cite{Shen2}. The Berwald-Riemann
curvature tensor $K^i_{\ jkl}$ is defined by \[ K^i_{\
jkl}=\frac{1}{3}\Big\{R^i_{\ k.l.j}-R^i_{\ l.k.j}\Big\}.
\]
The {\it Ricci scalar} is denoted by ${\bf Ric}$ it is defined by
${\bf Ric}:=K^k_{\ k}$. A Finsler space $(M,F)$ is called an {\it
Einstein space} if there is function $\lambda$ defined on $M$ such
that ${\bf Ric}=\lambda(x)F^2$. Notice that, the usual Ricci tensor
$K_{jl}$ has an antisymmetric part denoted by $\mathcal{R}$ and
defined by $\mathcal{R}_{jl}=\big(K_{\ jl}-K_{\ lj}\big)/2$. The
locally projectively flat Einstein Randers metric are locally
characterized in the following theorem, see \cite{ChenMoShen}.
\begin{thm}
\label{local characterization} Let $F=\alpha+\beta$ be a locally
projectively flat Randers metric on an n-dimensional manifold $M$.
Suppose that $F$ has constant Ricci curvature ${\bf Ric} = (n -
1)\lambda F^2$. Then $\lambda\leq0$. Further, if $\lambda = 0$, F is
locally Minkowskian. If $\lambda = -{1\over4}$, $F$ can be expressed
in the following form:
\[
F(x,y)=\frac{\sqrt{|y|^2-(|x|^2|y|^2-\langle
x,y\rangle^2)}}{1-|x|^2}\pm\frac{\langle
x,y\rangle}{1-|x|^2}\pm\frac{\langle a,y\rangle}{1+\langle
a,x\rangle},
\]
where, $y\in T_x\mathbb{R}^n,\ a\in\mathbb{R}^n, |a|<1$.
\end{thm}

\section{C-projective vector fields on Finsler spaces}
Every vector field $V$ on a manifold $M$ induces naturally an
infinitesimal coordinate transformations on $TM$ given by
$(x^i,y^i)\longrightarrow(\bar{x}^i,\bar{y}^i)$ is given loclly by
\begin{equation}
\bar{x}^i=x^i+V^idt,\ \ \ \ \ \bar{y}^i=y^i+y^k\frac{\partial
V^i}{\partial x^k}dt.
\end{equation}
This is in fact \textit{the complete lift} of $V$ (see for example
\cite{Yano}) to a vector field on $TM_0$ denoted by $\hat{V}$ and
given by
\begin{equation}
\hat{V}=V^i\frac{\partial}{\partial x^i}+y^k\frac{\partial
V^i}{\partial x^k}\frac{\partial}{\partial y^i}.
\end{equation}
Notice that, $\mathcal{L}_{\hat{V}}y^i=0$,
$\mathcal{L}_{\hat{V}}dx^i=0$ and the differential operators
$\mathcal{L}_{\hat{V}}$, $\frac{\partial}{\partial x^i}$, exterior
differential operator $d$ and $\frac{\partial}{\partial y^i}$
commute.

\bigskip

Now, let us suppose that $(M,F)$ be a Finsler space. A vector field
$V$ on $M$ is said to be \textit{projective}, if there is a function
$P$ (called the \textit{projective factor}) on $TM_0$ such that
$\mathcal{L}_{\hat{V}}G^i=Py^i$, see \cite{A1}. It is known that,
given any projective vector field $V$, its local flow $\{\phi_t\}$
associated to $V$ is a projective transformation, namely, $\phi_t$
sends forward geodesics to forward geodesics and vice-versa.  The
collection of the projective vector fields on a Finsler space
$(M,F)$ is denoted by $proj(M,F)$ and is classically known to be a
finite dimensional Lie algebra with respect to the usual Lie bracket
$[,]$. If $V$ is a projective vector field, then the following
identities are known  by \cite{A1, NajafiTayebi}:
\begin{eqnarray}
\mathcal{L}_{\hat{V}}G^i_{\ k}&=&P\delta^i_{\ k}+P_{k} y^i, \label{proj 3}\\
\mathcal{L}_{\hat{V}}G^i_{\ jk}&=&\delta^i_{\ j}P_k+\delta^i_{\
k}P_j+y^iP_{kj}, \label{proj 4}\\
2\mathcal{L}_{\hat{V}}{\bf E}_{jl}&=&(n+1)P_{jl}, \label{proj 6}\nonumber\\
2\mathcal{L}_{\hat{V}}{\bf H}_{jl}&=&(n+1)P_{jl|m}y^m, \label{proj 6 H}\nonumber\\
\mathcal{L}_{\hat{V}}K^i_{\ jkl}&=&\delta^i_{\
j}(P_{l|k}-P_{k|l})+\delta^i_{\ l}P_{j|k}-\delta^i_{\
k}P_{j|l}+y^i(P_{l|k}-P_{k|l})_{.j},\label{proj 10}\nonumber\\
\mathcal{L}_{\hat{V}}K_{\
jl}&=&P_{l|j}-nP_{j|l}+P_{lj|0},\label{proj 11}
\end{eqnarray}
where, $P_i=P_{.i}$ and $P_{ij}=P_{i.j}$, etc and $|$ denotes the
horizontal derivative with respect to the Berwald connection. A
projective vector field is said to be \textit{affine} if $P=0$. The
collections of the affine and Killing vector fields on a Finsler
space $(M,F)$ are denoted by $aff(M,F)$ and $k(M,F)$, respectively.
It is well-know that every Killing vector field is affine and every
affine vector field is projective. Thus, it is clear that
$k(M,F)\subseteq aff(M,F)\subseteq proj(M,F)$. Recall that, given
any projective vector field $V$ on a Riemannian spaces, the
projective factor $P=P(x,y)$ is linear with respect to $y$ and thus,
it is the natural lift of a 1-form on $M$ to a function on $TM_0$;
Notice that, by \eqref{proj 11}, the projective factor is actually a
closed 1-form for any Riemannian space (that is to say that
$P_{i|j}=P_{j|i}$), while these issue is a non-Riemannian feature in
a Finslerian media. Consider the following conventional definitions
of a projective vector field $V$; $V$ is said to be (cf.
\cite{RafieRad2011,RafieRad2013})
\begin{description}
   \item(i) {\it special} if $\mathcal{L}_{\hat{V}}{\bf E}=0$, or equivalently, $P(x,y)=P_i(x)y^i$.
   \item(ii) {\it C-projective} if $P_{i|j}=P_{j|i}$
   \item(iii) {\bf H}-{\it invariant} if $\mathcal{L}_{\hat{V}}{\bf H}=0$, equivalently, $P_{jk|l}=P_{jl|k}$
\end{description}
A projective vector field on a Riemannian manifold is simultaneously
a special and a C-projective. Every projective vector field on a
weakly Berwald space (i.e. ${\bf E}=0$ is special and every special
projective vector field on a Randers space of constant non-zero {\bf
S}-curvature is C-projective, cf. \cite{RafieRad2011,RafieRad2013}.
Now let us suppose that $\mathcal{R}_{lj}:=\big(K_{\ lj}-K_{\
jl}\big)/2$ denotes the anti-symmetric part of the usual Ricci
curvature $K_{jl}=K^i_{jil}$. Then by \eqref{proj 11}, the following
equation results immediately:
\begin{eqnarray}
\mathcal{L}_{\hat{V}}\mathcal{R}_{jl}&=&\Big(P_{l|j}-nP_{j|l}+P_{lj|0}-P_{j|l}+nP_{l|j}-P_{jl|0}\Big)/2\nonumber\\
&=&\frac{(n+1)}{2}\big(P_{l|j}-P_{j|l}\big).\label{C proj 12}
\end{eqnarray}
Using \eqref{C proj 12} we obtain the following characterization of
the C-projective vector fields:
\begin{cor}
Let us suppose that $(M,F=\alpha+\beta)$ be a Randers space and
$\mathcal{R}_{jl}=\big(K_{\ jl}-K_{\ lj}\big)/2$ denotes the usual
Ricci tensor. A projective vector field $V$ is C-projective if and
only if $\mathcal{L}_{\hat{V}}\mathcal{R}_{jl}=0$.
\end{cor}
The collection of the C-projective vector fields on a Finsler space
$(M,F)$ is denoted by $cproj(M,F)$. It is clear that
$cproj(M,F)\subseteq proj(M,F)$ and thus, $\dim\Big(cp(M,F)\Big)\leq
\dim\Big(p(M,F)\Big)\leq n(n+2)$. The case of maximum dimension
$n(n+2)$ is known on Randers paces by the following result:
\begin{thm}
\label{maximum} (Rafie-Rad and Rezaei, \cite{RafieRad-Rezaei2012})
 A Randers metric $F=\alpha+\beta$ on a manifold $M$ of dimension $n,\ (n\geq3)$ is projective if and only if $proj(M,F)$ has (locally)
dimension $n(n+2)$.
\end{thm}
The quotient $Q(M,F):=proj(M,F)/cproj(M,F)$ demonstrates a
non-Riemannian feature. The following result present shows that that
generalized Funk metrics on the  are examples whose projective
algebra consists only of the C-projective vector fields, namely
$Q(M,F):=proj(M,F)/cproj(M,F)=\{0\}$.

\bigskip

\begin{thm}
(Rafie-Rad, \cite{RafieRad2011}) Let $F=\alpha+\beta$ be an
$n$-dimensional Randers space of nonzero constant S-curvature. If
$F$ is projective (i.e. locally projectively flat) then, every
projective vector field is C-projective.
\end{thm}

\bigskip
\begin{rem}
\label{rem genfunk} The generalized Funk metric is a Randers metric
defined on the Euclidean unit ball $\mathbb{B}^n(1)$ as follows:
\begin{equation}
\label{local form} F(x,y)=\frac{\sqrt{|y|^2-(|x|^2|y|^2-\langle
x,y\rangle^2)}}{1-|x|^2}\pm\frac{\langle
x,y\rangle}{1-|x|^2}\pm\frac{\langle a,y\rangle}{1+\langle
a,x\rangle},
\end{equation}
where, $y\in T_x\mathbb{R}^n,\ a\in\mathbb{R}^n, |a|<1$.The above
metric is of nonzero constant S-curvature $S=\pm\frac{(n+1)}{2}F$
and theorem above ensures that they are explicit examples for which,
very projective vector field is C-projective. This fact also may
complete the half part of the proof of theorem \ref{mainthm2}.
\end{rem}
\bigskip
There are several projectively invariant tensors in Finsler geometry
such as the {\it Douglas} curvature $D=D^i_{\
jkl}\frac{\partial}{\partial x^i}\otimes dx^j\otimes dx^k\otimes
dx^l$ and \textit{Weyl} curvature $W=W^i_{\
jkl}\frac{\partial}{\partial x^i}\otimes dx^j\otimes dx^k\otimes
dx^l$ and the recently introduced quantity
$\widetilde{W}=\widetilde{W}{^i_{\ jkl}}\frac{\partial}{\partial
x^i}\otimes dx^j\otimes dx^k\otimes dx^l$ in \cite{NajafiTayebi}.
The tensors $D$, $W$ and $\widetilde{W}$ are defined as follow:
\begin{eqnarray}
D^i_{\ jkl}&=&\frac{\partial^3}{\partial y^j\partial y^k\partial y^l}\{G^i-\frac{1}{n+1}G^m_{\ m}y^i\},\label{Douglas}\nonumber\\
W^i_{\ jkl}&=&K^i_{\ jkl}-\frac{1}{n^2-1}\{\delta^i_{\
j}(\hat{K}_{kl}-\hat{K}_{lk})+(\delta^i_{\
k}\hat{K}_{jl}-\delta^i_{\
l}\hat{K}_{jk})+y^i(\hat{K}_{kl}-\hat{K}_{lk})_{.j}\},\nonumber\\
\widetilde{W}{^i_{\ jkl}}&=&K^i_{\ jkl}-\frac{\delta^i_l}{1-n^2}\Big\{\hat{K}_{jk}+\frac{n}{n+1}y^r(K_{jr.k}-K_{jk.r})\Big\}\nonumber\\
&
&+\frac{\delta^i_k}{1-n^2}\Big\{\hat{K}_{jl}+\frac{n}{n+1}y^r(K_{jr.l}-K_{jl.r})\Big\}\nonumber
\end{eqnarray}
where, $\hat{K}_{jk}=nK_{jk}+K_{kj}+y^rK_{kr.j}$. Notice that,
Akbar-Zadeh has already introduced a tensor which is just invariant
by a sub-group of projective transformations, not all of them
\cite{A1}. In fact, this is a non-Riemannian generalization of
Weyl's curvature. It is denoted by $\overset{*}{W}{^i_{\ jkl}}$ and
is defined by:
\[
\overset{*}{W}{^i_{\ jkl}}=K^i_{\ jkl}-\frac{1}{n^2-1}\{\delta^i_{\
k}(nK_{jl}+K_{lj})-\delta^i_{\ l}(nK_{jk}+K_{kj})+(n-1)\delta^i_{\
j}(K_{kl}-K_{lk})\}.
\]
The quantities $W$, $\widetilde{W}$ and $\overset{*}{W}$ are
invariant objects in the following invoice:
\begin{description}
  \item(i) $W$ is invariant under every projective vector field,
  \item(ii) $\widetilde{W}$ is invariant under every C-projective vector field,
  \item(iii) $\overset{*}{W}$ is invariant under every special projective vector field.
\end{description}

The projective vector field are variously characterized in many
contexts such as \cite{A1}. Some characterization of projective
vector fields in a Randers space $(M,F=\alpha+\beta)$ are give by in
terms of $\alpha$ and $\beta$ as follows:

\begin{thm}
\label{good} (Rafie-Rad and Rezaei, \cite{RafieRadRezaei2011})
\label{characterization} A vector field $V$ is projective on a
Randers space $(M,F=\alpha+\beta)$ if and only if $V$ is projective
in $(M,\alpha)$
 and $\mathcal{L}_{\hat{V}}(\alpha s^i_{\ j})=0$.
\end{thm}

There are also some characterization of special projective vector
fields:
\begin{thm}
\label{lem1} (Rafie-Rad and Rezaei,
\cite{RafieRad2011,RafieRad2013}) Let $(M,F=\alpha+\beta)$ be a
non-zero constant {\bf S}-curvature Randers space and $V\in\chi(M)$.
Then, $V$ is a special projective vector field if and only if $V$ is
a Killing vector field on $(M,\alpha)$ and
$\mathcal{L}_{V}s_{ij}=0$.
\end{thm}

\begin{rem}
\label{proj invariant} The quantity $\alpha s^i_{\
j}\frac{\partial}{\partial x^i}\otimes dx^j$ is a projectively
invariant tensor.
\end{rem}
The following result describes all special projective vector fields
on a Randers space  of isotropic {\bf S}-curvature; The case of
constant {\bf S}-curvature has already been proved in
\cite{RafieRad2011,RafieRad2013}.
\begin{thm}
\label{SP main1}
 Let $(M,F=\alpha +\beta )$ be a Randers space of
isotropic {\bf S}-curvature ${\bf S}=(n+1)c(x)F$. A vector field
$V\in \chi (M)$ is special projective if and only if there is a
1-form
$P=P_i(x)dx^i$ on $M$ such that:\\
(i) $\mathcal{L}_{\hat{V}}G^{i}_{\alpha}=\left( (V.c)\beta
+c(x)\mathcal{L}_{\hat{V}}\beta
+\mathcal{L}_{\hat{V}}s_{0}+P\right)y^{i}$,\\
(ii) $(2(V.c))\alpha ^{2}+c(x)t_{00}=0.$
\end{thm}
\textit{Proof.} Let $(M,F=\alpha +\beta )$ is a Randers space of
isotropic curvature ${\bf S}=(n+1)c(x)F$. Following Theorem
\ref{characterization}, $V$ is a special projective vector field on
$M$ iff we have $\mathcal{L}_{\hat{V}}G^{i}=P(x,y)y^{i}$ and
$\mathcal{L}_{\hat{V}}(\alpha s_{0}^{i})=0$, where, $P$ is a 1-form
on $M$. The geodesic spray coefficients $G^i$ are given by
\eqref{Gi}; Hence, it follows:
\begin{equation}
\label{Sh} \mathcal{L}_{\hat{V}}\left( G^{i}_{\alpha}+(c(x)(\alpha
-\beta )-s_{0})y^{i}\right)=Py^i
\end{equation}
Let us assume $t_{ij}:=\mathcal{L}_{\hat{V}}a_{ij}$ and
$V.c:=\mathcal{L}_{\hat{V}}c$; Thus, $\alpha
^{2}=a_{ij}(x)y^{i}y^{j}$ and $t_{00}=\mathcal{L}_{\hat{V}}\alpha
^{2}$ and then $\mathcal{L}_{\hat{V}}\alpha
=\frac{t_{00}}{2\alpha}$. Now \eqref{Sh} is equivalent to the
following equation:
\begin{equation}
\label{W} \mathcal{L}_{\hat{V}}G^{i}_{\alpha}+\left((V.c)(\alpha
-\beta )+c(x)\dfrac{t_{00}}{2\alpha}-c(x)\mathcal{L}_{\hat{V}}\beta
-\mathcal{L}_{\hat{V}}s_{0}-P\right)y^{i}=0
\end{equation}
It is clear that that terms $(\alpha^2-\beta^2)$ and
$c(x)\mathcal{L}_{\hat{V}}\beta$ and $\mathcal{L}_{\hat{V}}s_{0}$
are polynomial with degree 2, 1 and 1, respectively. Now,
multiplying the two side of \eqref{W} by $2\alpha$ and simplifying
the terms, it results $Rat^i+\alpha Irrat^i=0,\ (i=1,...,n)$, where,
\begin{eqnarray*}
Irrat ^{i}&=&2\mathcal{L}_{\hat{V}}G^{i}_{\alpha}-\Big(2(V.c)\beta
+2c(x)\mathcal{L}_{\hat{V}}\beta
+2\mathcal{L}_{\hat{V}}s_{0}+2P\Big)y^i,\\
Rat^{i}&=&\left(2(V.c)\alpha ^{2}+c(x)t_{00}\right)y^{i}.
\end{eqnarray*}
Therefore, $V$ is a special projective vector field if and only if
$Rat^i=0$ and $Irrat^i=0$ and theorem is proved. $\Box$\\

\bigskip

Let us suppose that $V$ is a projective vector field on a Randers
space $(M,F=\alpha+\beta$. Therefore,
$\mathcal{L}_{\hat{V}}G^i=Py^i$ and
$\mathcal{L}_{\hat{V}}G^i_\alpha=\eta y^i$, where, $P$ and $\eta$
denote the projective factors for $F$ and $\alpha$, respectively.
Now, by \eqref{GiS} and Remark \ref{proj invariant}, we have
\begin{eqnarray*}
\mathcal{L}_{\hat{V}}G^i&=&Py^i\\
&=&\mathcal{L}_{\hat{V}}\Big\{G^i_\alpha+\Big(\frac{\bf
S}{n+1}+\rho_0\Big)y^i+\alpha s^i_{\ 0}\Big\}\\
&=&\Big\{\eta+\mathcal{L}_{\hat{V}}\Big(\frac{\bf
S}{n+1}+\rho_0\Big)\Big\}y^i=Py^i
\end{eqnarray*}
\begin{rem}
\label{Proj factor} The projective factor for every projective
vector field $V$ on a Randers space $(M,F=\alpha+\beta)$ is given by
$P=\eta+\mathcal{L}_{\hat{V}}\Big(\frac{\bf S}{n+1}+\rho_0\Big)$,
where, $P$ and $\eta$ denote the projective factors for $F$ and
$\alpha$, respectively.
\end{rem}
\bigskip

\section{Proof of main Theorems.}
\bigskip
\textbf{Proof of Theorem \ref{mainthm1}:} \textit{Proof of
implication $(1)\Rightarrow (2)$:} Let us suppose that $V$ is a
C-projective vector field, i.e. $\mathcal{L}_{\hat{V}}G^i=Py^i$ and
$P_{i|j}=P_{j|i}$. By Remark \ref{Proj factor}, the projective
factor $P$ is given by $P=\eta+\mathcal{L}_{\hat{V}}\Big(\frac{\bf
S}{n+1}+\rho_0\Big)$, where, $\eta$ is the projective factor for
$\alpha$; Recall that $\eta$ is a closed form. Hence, it results
that
\begin{eqnarray}\label{S 001}
P_{i|j}&=&\eta_{i|j}+\Big\{\mathcal{L}_{\hat{V}}\Big(\frac{{\bf
S}_{.i}}{n+1}+\rho_i\Big)\}_{|j}\\
&=&\eta_{i|j}+\frac{1}{n+1}\Big(\mathcal{L}_{\hat{V}}{\bf
S}_{.i}\Big)_{|j}+\Big(V.\rho\Big)_{;i|j}\nonumber
\end{eqnarray}
On the other hand, following \eqref{Cov der}, \eqref{proj 3} and
\eqref{proj 4}, the term $\mathcal{L}_{\hat{V}}{\bf S}_{.i|j}$ can
be obtained as follows:
\begin{eqnarray*}
\mathcal{L}_{\hat{V}}{\bf
S}_{.i|j}&=&\mathcal{L}_{\hat{V}}\Big\{\frac{\partial}{\partial
x^j}{\bf S}_{.i}-{\bf S}_{.i.r}G^r_{\ j}-{\bf
S}_{.r}G^r_{ij}\Big\}\\
&=&\frac{\partial}{\partial x^j}\mathcal{L}_{\hat{V}}{\bf
S}_{.i}-\big(\mathcal{L}_{\hat{V}}{\bf S}_{.i.r}\big)G^r_{\ j}-{\bf
S}_{.i.r}\mathcal{L}_{\hat{V}}G^r_{\
j}-\big(\mathcal{L}_{\hat{V}}{\bf S}_{.r}\big)G^r_{ij}-{\bf
S}_{.r}\mathcal{L}_{\hat{V}}G^r_{ij}\nonumber\\
&=&\Big(\mathcal{L}_{\hat{V}}{\bf S}_{.i}\Big)_{|j}-{\bf
S}_{.i.r}\big(P\delta^r_j+P_{j}y^r\big) -{\bf
S}_{.r}\big(P_{ij}y^r+P_{i}\delta^r_{\ j}+P_{j}\delta^r_{\
i}\big)\nonumber\\
&=&\Big(\mathcal{L}_{\hat{V}}{\bf S}_{.i}\Big)_{|j}-P{\bf
S}_{.i.j}-{\bf S}P_{ij}-{\bf S}_{.i}P_{j}-{\bf S}_{.j}P_{i}\nonumber
\end{eqnarray*}
and finally we obtain $\Big(\mathcal{L}_{\hat{V}}{\bf
S}_{.i}\Big)_{|j}$ as follows:
\begin{eqnarray}\label{S 002}
\Big(\mathcal{L}_{\hat{V}}{\bf
S}_{.i}\Big)_{|j}=\mathcal{L}_{\hat{V}}{\bf S}_{.i|j}+P{\bf
S}_{.i.j}+{\bf S}P_{ij}+{\bf S}_{.i}P_{j}+{\bf S}_{.j}P_{i}
\end{eqnarray}
By closeness of $\eta$ (i.e. $\eta_{i|j}=\eta_{j|i}$) and
$\Big(V.\rho\Big)_{;i|j}=\Big(V.\rho\Big)_{;j|i}$ and also plugging
$\Big(\mathcal{L}_{\hat{V}}{\bf S}_{.i}\Big)_{|j}$ from \eqref{S
002} in \eqref{S 001}, it results:
\begin{eqnarray}
P_{i|j}-P_{j|i}&=&\mathcal{L}_{\hat{V}}\Big\{\frac{1}{n+1}\big({\bf
S}_{.i|j}-{\bf
S}_{.j|i}\big)\Big\}=\mathcal{L}_{\hat{V}}\Sigma_{ij}\nonumber
\end{eqnarray}

\textit{Proof of the implication $(2)\Rightarrow (3)$:} Let us
suppose that $\mathcal{L}_{\hat{V}}\Sigma_{ij}=0$. Derivation with
respect to $y^k$ commutes with the Lie derivative operator; This
yields $\mathcal{L}_{\hat{V}}\Sigma_{ij.k}=0$. Due to \eqref{R3}, it
follows
\[
-\frac{2}{n+1}\mathcal{L}_{\hat{V}}H_{jk}=\mathcal{L}_{\hat{V}}\big(y^i\Sigma_{ij.k}\big)=y^i\mathcal{L}_{\hat{V}}\Sigma_{ij.k}=0.
\]
Taking into account $\mathcal{L}_{\hat{V}}H_{jk}=0$ and using
\eqref{R4}, we obtain $\mathcal{L}_{\hat{V}}\Xi_{j.k}=0$. Finally,
by \eqref{R02} it results
\[
-\frac{1}{n+1}\mathcal{L}_{\hat{V}}\Xi_k=\mathcal{L}_{\hat{V}}\big(y^j\Xi_{j.k}\big)=-y^j\mathcal{L}_{\hat{V}}\Xi_{j.k}=0.
\]
\textit{Proof of the implication $(3)\Rightarrow (1)$:} Let us
suppose that $\mathcal{L}_{\hat{V}}G^i=Py^i$ and
$\mathcal{L}_{\hat{V}}\Xi_i=0$. By Remark \ref{Proj factor},
$P=\eta+\mathcal{L}_{\hat{V}}\Big(\frac{\bf S}{n+1}+\rho_0\Big)$,
where, $\eta$ is the projective factor for $\alpha$ and is known to
be closed. Now by \eqref{R02}, \eqref{R5}, \eqref{S 001} and
\eqref{S 002} it follows that
\begin{eqnarray*}
y^mP_{i|m}-P_{|i}&=&y^m\Big(P_{i|m}-P_{m|i}\Big)\\
&=&-y^m\mathcal{L}_{\hat{V}}\Sigma_{mi}\\
&=&-\mathcal{L}_{\hat{V}}\big(y^m\Sigma_{mi}\big)\\
&=&-\frac{1}{n+1}\mathcal{L}_{\hat{V}}\big(y^m\Xi_{mi}\big)\\
&=&\frac{1}{n+1}\mathcal{L}_{\hat{V}}\Xi_i=0
\end{eqnarray*}
Recall Thus, $y^mP_{i|m}-P_{|i}=0$ and a derivation of the two sides
with respect to $y^j$ yields $P_{i|j}-P_{ij|m}y^m-P_{j|i}=0$. It
follows that $P_{i|j}-P_{j|i}=P_{ij|m}y^m$. The left hand is
anti-symmetric while the right hand is symmetric with respect to
indices $i$ and $j$. Thus, $P_{i|j}-P_{j|i}=0$ and $V$ is a
C-projective vector field. $\Box$
\bigskip

\begin{rem} \label{rem final} Suppose that $V$ is an arbitrary projective vector fields on the
Finsler space $(M,F)$. Taking into account the equation below
\begin{eqnarray}
P_{i|j}-P_{j|i}&=&\mathcal{L}_{\hat{V}}\Sigma_{ij},\nonumber
\end{eqnarray}
appearing in the proof of the implication $(1)\Rightarrow (2)$ and
using \eqref{proj 11}, we obtian
\[
\mathcal{L}_{\hat{V}}\Big\{\mathcal{R}_{jl}-\frac{(n+1)}{2}\Sigma_{jl}\Big\}=0.
\]
Thus, the tensor
$Z_{jl}:=\mathcal{R}_{jl}-\frac{(n+1)}{2}\Sigma_{jl}$ is a
projective invariant.
\end{rem}

\textbf{Proof of Theorem \ref{mainthm2}:} Suppose that the
C-projective algebra $cp(M,F)$ has maximum dimension $n(n+2)$. It
follows immediately that the whole projective algebra $p(M,F)$ has
dimension $n(n+2)$ and thus $p(M,F)=cp(M,F)$. By Theorem
\ref{maximum}, it follows that $F$ is a projective metric;
Equivalently, $\alpha$ has constant sectional curvature and
$s_{ij}=0$. Notice that, in this case $F$ and $\alpha$ are
projectively equivalent. Now by Theorem \ref{mainthm1}, given any
projective vector field $V$-which is now a C-projective vector field
too-we have $\mathcal{L}_{\hat{V}}\Sigma_{ij}=0$. In particular,
every Killing vector field for $\alpha$ is also a projective for $F$
(and also is C-projective). $\alpha$ has constant sectional
curvature say $k$, hence, the algebra of Killing vectors of
$(M,\alpha)$ has maximal dimension $n(n+1)/2$. It is also well-known
that the Killing vector field $V$ is locally of the form
\begin{equation}
\label{lie1} V^i=Q^i_{\ k}x^k+C^i+k\langle x,C\rangle x^i,
\end{equation}
where, $C$ is an arbitrary constant vector and $Q^i_{\ k}$ is an
arbitrary constant skew-symmetry bilinear form. On the other hand,
$\mathcal{L}_{\hat{V}}\Sigma_{ij}=0$ gives
\begin{eqnarray}
\label{lie2} \mathcal{L}_{\hat{V}}\Sigma_{ij}&=&\frac{\partial
V^k}{\partial x^j}\Sigma_{ik}+\frac{\partial V^k}{\partial
x^i}\Sigma_{kj}+V^k\Sigma_{ij,k}+y^m\frac{\partial V^k}{\partial
x^m}\Sigma_{ij.k}=0;
\end{eqnarray}
where, the subscript $_{,k}$ denotes the derivative with respect to
$\frac{\delta}{\delta x^k}$. From (\ref{lie1}), we obtain
\begin{eqnarray}
\frac{\partial V^k}{\partial x^j}&=&Q^k_{\ j}+k\langle x,C\rangle\delta^k_{\ j}+kC^jx^k,\label{lie3}\\
\frac{\partial V^k}{\partial x^i}&=&Q^k_{\ i}+k\langle
x,C\rangle\delta^k_{\ i}+kC^ix^k\label{lie4}.
\end{eqnarray}
Plugging the terms $\frac{\partial V^k}{\partial x^j}$ and
$\frac{\partial V^k}{\partial x^i}$ from (\ref{lie3}) and
(\ref{lie4}) in (\ref{lie2}) and taking into account $C=0$, it
results:
\begin{eqnarray} Q^k_{\ j}\Sigma_{ik}+Q^k_{\ i}\Sigma_{kj}+Q^k_{\
l}T^l_{\ ijk}=0, \label{lie6}
\end{eqnarray}
where, $Q=\Big(Q^k_{\ j}\Big)$ is an arbitrary skew-symmetric matrix
and \[ T^l_{\ ijk}=x^l\Sigma_{ij,k}+y^l\Sigma_{ij.k}.\] Consider two
fixed distinct indices $l_0$ and $k_0$ such that $Q^{k_0}_{\
l_0}=-Q^{l_0}_{\ k_0}=1$ and $Q^k_{\ l}=0$ if $k\neq k_0$ or $l\neq
l_0$. Given any indices $i$ and $j$ such that $i,j\neq l_0$, we have
\[
Q^k_{\ j}\Sigma_{ik}=0,\ \ \ Q^k_{\ i}\Sigma_{kj}=0,\ \ \ Q^k_{\
l}x^l=\left\{
                                                                  \begin{array}{ll}
                                                                  x^{l_0}, \ \  k=k_0\\
                                                                  -x^{k_0}, \ \ k=l_0 \\
                                                                   0,\ \ \ \textrm{otherwise}.
                                                                  \end{array}
                                                                  \right.
                                                                  ,\ \ \ Q^k_{\ l}y^l=\left\{
                                                                  \begin{array}{ll}
                                                                  y^{l_0}, \ \  k=k_0\\
                                                                  -y^{k_0}, \ \ k=l_0 \\
                                                                   0,\ \ \ \textrm{otherwise}.
                                                                  \end{array}
                                                                  \right.
\]
the equation (\ref{lie6}) becomes
\begin{equation}
\label{L1} T^{l_0}_{\ ijk_0}-T^{k_0}_{\ ijl_0}=0.
\end{equation}
It follows that, \eqref{L1} holds if $i,j\neq k$. Now, Fix two
distinct indices $i$ and $j$ and consider the matrix $Q$ given by
$Q^i_{\ j}=-Q^j_{\ i}=1$ and $Q^k_{\ l}=0$ if $k\neq i$ or $l\neq
l_0$. Observe that for the matrix $Q$ we have
\[
Q^k_{\ j}\Sigma_{ik}=Q^i_{\ j}\Sigma_{ii}=0,\ \ \ Q^k_{\
i}\Sigma_{kj}=Q^j_{\ i}\Sigma_{jj}=0,
\]
\[ Q^k_{\ l}x^l=\left\{
                                                                  \begin{array}{ll}
                                                                  x^{j}, \ \  k=i\\
                                                                  -x^{i}, \ \ k=j \\
                                                                   0,\ \ \ \textrm{otherwise}.
                                                                  \end{array}
                                                                  \right.
                                                                  ,\ \ \ Q^k_{\ l}y^l=\left\{
                                                                  \begin{array}{ll}
                                                                  y^{j}, \ \  k=i\\
                                                                  -y^{i}, \ \ k=j \\
                                                                   0,\ \ \ \textrm{otherwise}.
                                                                  \end{array}
                                                                  \right.
\]
and the equation (\ref{lie6}) becomes
\begin{equation}
\label{L2} T^j_{\ iji}-T^i_{\ ijj}=0
\end{equation}
Therefore, from \eqref{L1} and \eqref{L2}, it follows that given any
skew-symmetric matrix $Q=\big(Q^k_{\ l}\big)$, we have $T^l_{\
ijk}=T^k_{\ ijl}$ and thus,
\begin{equation}
\label{L3} Q^k_{\ l}T^l_{\ ijk}=0.
\end{equation}
Plugging \eqref{L3} in \eqref{lie6}, we obtain
\begin{eqnarray} Q^k_{\ j}\Sigma_{ik}+Q^k_{\ i}\Sigma_{kj}=0, \label{lie7}
\end{eqnarray}

Now, let $i\neq j$ and $k_0\neq i,j$ and $Q^{k_0}_{\ i}=-Q^i_{\
k_0}=1$ and $Q^k_{\ l}=0$ if $k\neq k_0$ or $l\neq i$. Thus,
equation (\ref{lie7}) can be written as follows: $Q^k_{\
j}\Sigma_{ik}+Q^k_{\ i}\Sigma_{kj}=\Sigma_{k_0j}=0$. Since $i$ and
$j$ are arbitrarily chosen, hence $\Sigma_{ij}=0$. Moreover,
$\Xi_i=-y^j\Sigma_{ij}=0$. By \eqref{R3}, it results that
$H_{jk}=-y^i\Sigma_{ij.k}=0$ and the ${\bf H}$-curvature vanishes.
We summarize the results as $F$ is a projective Randers metric with
constant flag curvature ${\bf K}=\lambda$ and subsequently, it is of
constant Ricci curvature ${\bf Ric}=(n-1)\lambda F^2$. Now, by
Theorem \ref{local characterization}, $F$ is locally Minkowskian or
up to a re-scaling, $F$ is locally isometric to the generalized Funk
metric given in \eqref{local form}.\\

 Conversely, the generalized Funk
$F$ metrics on the Euclidean unit ball $\mathbb{B}^n(1)$ have
constant S-curvature ${\bf S}=\pm\frac{(n+1)}{2}F$. This follows
immediately that given any projective vector field $V$ on
$\mathbb{B}^n(1)$, we have
\begin{eqnarray*}
\Sigma_{ij}&=&\frac{1}{n+1}\Big\{{\bf S}_{.i|j}-{\bf S}_{.j|i}\Big\}\\
&=&\pm\frac{1}{2}\Big\{F_{.i|j}-F_{.j|i}\Big\}=0.
\end{eqnarray*}
Now by Remark \ref{rem final} it follows that $P_{i|j}=P_{j|i}$ and
$V$ is a C-projective vector field; Thus, $proj(M,F)=cproj(M,F)$.
The generalized Funk metrics are projective and it is known that
$proj(M,F)=proj(M,\alpha)$. $\alpha$ is of constant sectional
curvature and in this case, $proj(M,\alpha)$ has local dimension
$n(n+2)$ and so does $cproj(M,F)=proj(M,F)$. $\ \ \ \ \ \Box$

\end{document}